\documentclass[10pt,A4,conference]{IEEEtran}
\usepackage{geometry}
 \usepackage{authblk}
\usepackage{tikz}  
\usetikzlibrary{calc,intersections}
 \usepackage{verbatim}
\usepackage{textcomp}
\usepackage[T1]{fontenc}
\usepackage{bbm}
\usepackage{multirow}
\usepackage{dsfont}
\usepackage{cite}
\usepackage{epsfig}
\usepackage{float}
\usepackage{enumerate}
\usepackage{balance}
\usepackage{color}
\usepackage{subcaption}
\usepackage{caption}
\usepackage{algorithm}
\usepackage[noend]{algpseudocode}
\algrenewcommand\Return{\State \algorithmicreturn{} } 
\usepackage{hyperref}
\hypersetup{
    colorlinks=true,
    linkcolor=blue,
    filecolor=magenta,      
    urlcolor=cyan,
		citecolor=blue,
}
\usepackage{amssymb}
\usepackage{amsthm}
\usepackage{array}
\usepackage{graphicx}
\usepackage{lettrine}
\usepackage{epsf} 
\usepackage{psfrag}
\usepackage[usenames,dvipsnames]{pstricks}
\usepackage{pst-grad} 
\usepackage{pst-plot} 
\usepackage[space]{grffile} 
\usepackage{etoolbox} 
\usepackage{url} 

\usepackage{multirow}

\usepackage [ short ]{ optidef }
\usepackage[capitalize]{cleveref}
\usepackage[colorinlistoftodos,bordercolor=orange,backgroundcolor=orange!20,linecolor=orange,textsize=scriptsize]{todonotes}
\Crefname{figure}{\Fig.}{Figs.}

\begin{document}

\title{Age-of-Information in UAV-assisted Networks: a Decentralized Multi-Agent Optimization}

 
\author[]{Mouhamed Naby Ndiaye}
\author[]{El Houcine Bergou}
\author[]{Hajar El Hammouti}
\affil[]{College of Computing, Mohammed VI Polytechnic University (UM6P), Benguerir, Morocco, \authorcr
emails: {\{naby.ndiaye,elhoucine.bergou,hajar.elhammouti\}@um6p.ma}}


\maketitle

\begin{abstract}
Unmanned aerial vehicles (UAVs) are a highly promising technology with diverse applications in wireless networks. One of their primary uses is the collection of time-sensitive data from Internet of Things (IoT) devices. In UAV-assisted networks, the Age-of-Information (AoI) serves as a fundamental metric for quantifying data timeliness and freshness. In this work, we are interested in a generalized AoI formulation, where each packet's age is weighted based on its generation time. Our objective is to find the optimal UAVs' trajectories and the subsets of selected devices such that the weighted AoI is minimized. To address this challenge, we formulate the problem as a Mixed-Integer Nonlinear Programming (MINLP), incorporating time and quality of service constraints. To efficiently tackle this complex problem and minimize communication overhead among UAVs, we propose a distributed approach. This approach enables drones to make independent decisions based on locally acquired data. Specifically, we reformulate our problem such that our objective function is easily decomposed into individual rewards. The reformulated problem is solved using a distributed implementation of Multi-Agent Reinforcement Learning (MARL). Our empirical results show that the proposed decentralized approach achieves results that are nearly equivalent to a centralized implementation with a notable reduction in communication overhead.

\end{abstract}
\IEEEoverridecommandlockouts
\IEEEpeerreviewmaketitle


\section{Introduction}
Unmanned aerial vehicles (UAVs), also known as drones, are used in wireless communications for a variety of applications such as collecting data and providing on-demand connectivity. However, their deployment in a multi-agent environment poses several challenges such as resource allocation, trajectory design, and cooperation between drones~\cite{bajracharya20226g,hammouti2018air}. In this paper, we focus on the use of UAVs to periodically collect time-sensitive data from IoT devices and transmit it to a server for analysis and decision-making. To measure the freshness of the collected data, we are interested in the concept of Age-of-Information (AoI)~\cite{kosta2017age,10333749} which quantifies the time elapsed since the last data update was collected. Our main objective is to minimize the weighted AoI of the network over a period of time. Since UAVs operate in a common environment, we need to carefully design their behaviors and interactions to achieve the common goal of minimizing the weighted AoI.

\subsection{Related Work}

Optimizing a common objective in multi-UAV networks is a challenging task for two main reasons. First, it requires finely tuned UAVs actions while respecting resource and energy constraints. Second, it entails the complex coordination between UAVs which should involve a minimum communication overhead. To address these challenges, Multi-Agent Reinforcement Learning (MARL), an extension of Reinforcement Learning (RL) to the multi-agent setting, has been proposed~\cite{arulkumaran2017deep,Cui2020MultiAgentRL,Seid2021MultiAgentDF,9900429}. In~\cite{Cui2020MultiAgentRL}, authors study the problem of dynamic resource allocation in multi-UAV networks. The studied problem is formulated as a stochastic game. To solve the optimization problem, the authors propose an RL-based approach where all agents share a common structure based on Q-learning. Similarly, in~\cite{Seid2021MultiAgentDF}, the authors study the deployment of a clustered multi-drone system to provide computational task offloading to IoT devices. They formulate their problem with the goal of minimizing the long-term computational cost in terms of delay and energy. Then, they propose a centralized method based on multi-agent deep reinforcement learning to minimize the overall computational cost of the network. In~\cite{9900429}, the authors propose a cooperative multi-agent algorithm with an actor-critic architecture to solve the problem of latency. The proposed critic architecture is shared between the agents.

In the previously cited works, MARL is implemented using a Centralized Training and Decentralized Execution (CTDE) strategy. In CTDE, the training and execution phases of the multi-agent system are handled differently. First, a centralized entity uses a global knowledge of the system to learn a coordinated policy for agents. Second, after the training, each agent operates independently and makes decisions based on its local observations. While CTDE achieves promising results, it requires the exchange of a large amount of information with the central entity, which may not be feasible in scenarios with limited communication bandwidth and latency. CTDE also does not scale well and becomes less efficient as the number of agents in the system increases.

On the other hand, a handful of works applying Decentralized Training and Decentralized Execution (DTDE) to MARL can be found~\cite{Lin2021DecentralizedPD,Hu2021DistributedMM,Hu2022DistributedAD}. In fact, in DTDE, both the training and execution phases are conducted in a decentralized manner. Each agent learns its policy independently and makes decisions without direct access to the global state or actions of other agents. For example, in \cite{Lin2021DecentralizedPD}, a decentralized trajectory planning is proposed for collision and obstacle avoidance in UAV networks.  The trajectory design of UAVs is also investigated in~\cite{Hu2021DistributedMM}, where energy-constrained drones serve dynamic users. The authors propose a distributed Value Decomposition-based Reinforcement Learning (VDRL) solution to find the trajectories of UAVs which maximize the fraction of served users. Similarly, a value decomposition RL is used in \cite{Hu2022DistributedAD} to address the problem of data pre-storing and routing in resource-constrained cube satellite networks.

None of the previously cited works studied the problem of AoI minimization. Furthermore, these works employ classical off-policy and on-policy algorithms to train agents which sometimes turn out to be unstable and data inefficient~\cite{Schulman2017ProximalPO}. To overcome this shortcoming, Proximal Policy Optimization (PPO) was introduced as an on-policy algorithm that benefits from both stability and data efficiency~\cite{Schulman2017ProximalPO,Yu2021TheSE}. In fact, PPO uses a special objective function that allows the policy to change rapidly, while still ensuring that it remains close to the old policy. This results in a more stable and data-efficient algorithm, which has been shown to be competitive with some of the widely used on- and off-policy algorithms.

The closest research to this work is our previously published paper~\cite{naby2023muti}, wherein we introduced a centralized Multi-Agent PPO (MAPPO) approach for minimizing the AoI. In contrast to our earlier work, this paper presents a decentralized implementation of MAPPO while relying on an approximated formulation of the studied optimization problem. Moreover, it relies on a more general formulation of the AoI where
a weight is attributed to the age of each packet according to its
generation time. Furthermore, a comparative analysis between the proposed decentralized approach and the centralized approach in~\cite{naby2023muti} is conducted.

\subsection{Contribution}
Our work focuses on a multi-UAV network in which a set of drones is deployed to collect time-sensitive data updates from IoT devices. Our goal is to optimally design the UAVs' trajectories and subsets of visited IoT devices, with the objective of minimizing the weighted AoI. Our contributions can be summarized as follows.

\begin{itemize}
\item We provide a general formulation of AoI where weight is attributed to the age of each packet according to its generation time. Each packet is assigned a weight, allowing for a greater emphasis on recently generated packets without ignoring previously generated ones. The problem is then formulated as
a Mixed-Integer Nonlinear Programming (MINLP), incorporating
time and quality of service constraints.

\item To solve the studied problem, we propose an approximation of the formulated problem that enables the decomposition of the reward, thus allowing the UAVs to learn optimal policies independently.
\item We propose a decentralized approach to training agents based on DTDE and PPO. Our approach aims to address the challenge of learning optimal policies for UAVs without the need for sharing information such as state, action, and reward data among the agents.
\item Finally, based on our simulation results, we show that the proposed decentralized approach achieves results that are nearly comparable to a centralized implementation, while significantly reducing communication overhead.
\end{itemize}

\subsection{Organization}
The remainder of the paper is organized as follows. In Section~\ref{Sys}, we describe the studied system model and introduce the AoI metric. The mathematical formulation of the problem is given in Section~\ref{Prob}. In Section~\ref{algo}, we describe in detail the proposed distributed MARL approach. Next, in Section~\ref{simu}, we assess the performance of the proposed approach. Finally, concluding remarks are provided in Section~\ref{Conc}.

\section{System Model}\label{Sys}
Consider a set of $I$ IoT devices, denoted by $\mathcal{I}$, which generate data updates periodically over a time span $T$. We assume that $T$ is partitioned into equal intervals of duration $\tau$. Within each interval, device $i$ generates new data at a period $k_i\tau$, where $k_i$ can take values from the set $\{1, 2, \ldots, K\}$, and $K=\frac{T}{\tau}$. We also assume that the data generated by each device is accumulated in a buffer until its collection by a UAV. 

\subsection{Age-of-Information Metric}
To capture the heterogeneous generation frequencies of IoT devices and accurately reflect the freshness of the collected information, we introduce AoI metric. The AoI quantifies the time elapsed between data generation and data collection for each IoT device. Let $A_{i}^n[t]$ represent the AoI of data at time interval $t$ generated by device $i$ during its $n^{\text{th}}$ period of generation, which corresponds to the $nk_i^{\text{th}}$ time interval. The AoI value $A_{i}^n[t]$ at time interval $t\geq 0$ can be calculated recursively using the following expression

\begin{equation*}\label{equa1}
    A_{i}^n[t]=\Bigg\{
    \begin{array}{ll}
   A_{i}^{n}[t-1]+\tau & \text{if} \sum\limits_{l=nk_i}^{t}\sum\limits_{u \in \mathcal{U}}\alpha_{iu}[t]=0 \text{ \& } nk_i\leq t\\
    0 & \text{otherwise},
    \end{array}
\end{equation*}
where $\mathcal{U}$ represents the set of UAVs, and $\alpha_{iu}[t]$ denotes a binary variable with $\alpha_{iu}[t]=1$ if UAV $u$ collects data from device $i$ during time interval $t$, and $\alpha_{iu}[t]=0$ otherwise, and $A_i^n[0]=0$. The AoI metric is formulated as follows: when data updates from device $i$ are not collected during time interval $t$ (i.e., $\sum \limits_{u \in \mathcal{U}}\alpha_{iu}[t]=0$), the AoI is increased by $\tau$. Conversely, when updates are collected (i.e., $\sum \limits_{u \in \mathcal{U}} \alpha_{iu}[t]=1$) or have not been generated yet (i.e., $nk_i \geq t$), the AoI is set to zero~\cite{ndiaye2022age,AbdElmagid2019DeepRL}.  Unlike existing works which focus on the AoI of the most recently generated packet~\cite{kadota2021age,bedewy2019age}, we adopt a more general formulation by considering all the packets in the buffer. In fact, in some applications, historical data may have an important impact on the process of prediction and decision-making, necessitating the evaluation of both older and recent data. To this end, we introduce a weight factor denoted as $w^n[t]$, which corresponds to the weight during time interval $t$ of a packet generated at period $n$. In our work, we assume that the weight factor increases as the generation time of the packet increases. This allows for a higher emphasis on the recently
generated packets while still considering previously generated
ones. Consequently, a packet generated at time $n$ holds less importance than another one generated at time $n+1$, i.e., $w^n[t]<w^{n+1}[t]$. A typical function of $w^n[t]$ is given by $w^n[t]=\gamma^{t-n}$, where $\gamma \in [0,1]$.

Finally, the weighted AoI of all the data generated by IoT device $i$ during the considered time, $f_i(\boldsymbol{\alpha}_i)$, is given by $f_i(\boldsymbol{\alpha}_i)= \sum \limits_{t=1}^K \sum\limits_{n=0}^{\lfloor \frac{K}{k_i} \rfloor}w^n[t]A_{i}^n[t]$, with $\boldsymbol{\alpha}_i=(\alpha_{iu}[t])_{\substack{u\in \mathcal{U} \\ t \in \{0,\dots,K-1\}}}$ a $U\times K$ matrix.


\subsection{Communication Model}
 Each IoT device communicates with a UAV through the air-to-ground channel. In our model, we adopt a block Rician-fading approach, where the propagation channels between any device-UAV pair remain constant during a time interval of length equal to or less than $\tau$. As a result, the air-to-ground channel between device $i$ and UAV $u$ is $h_{iu}[t]=d_{iu}^{-1}[t]\left(\sqrt{\frac{\Phi}{\Phi+1}} {\xi}^{\rm LoS}_{iu}[t]+\sqrt{\frac{1}{\Phi+1}} {\xi^{\rm NLoS}}_{iu}[t]\right)$, where $\Phi$ is the Rician factor, ${\xi}^{LoS}_{iu}[t]$ is the line-of-sight (LoS) component with $\left|{\xi}^{\rm LoS}_{iu}[t]\right|=1$, and ${\xi}^{\rm NLoS}_{iu}[t]$ the random non-line-of-sight (NLoS) component, which follows a Rayleigh distribution with mean zero and variance one. Finally, $d_{iu}[t]$ is the distance between the pair device-UAV during time interval $t$ that is given by $d_{iu}[t]=\sqrt{(x_{u}[t]-x_{i})^2+(y_{u}[t]-y_{i})^2+(H_{u})^2}$, where $(x_u[t],y_u[t],H_u)$ and $(x_i,y_i,0)$ are the 3D positions of UAV $u$ and device $i$ respectively. To avoid collisions, each UAV $u$ flies at a constant and different altitude  $H_u$.

We assume that devices use orthogonal frequency division multiple access (OFDMA) for their transmissions. Accordingly, the signal-to-noise ratio (SNR) of IoT device $i$ with respect to UAV $u$ is given by $\Gamma_{iu}[t]= P_{i}[t]\left|h_{iu}[t]\right|^{2} / \sigma^{2}$, where $\sigma^{2}$ is the variance of an additive white Gaussian noise and $P_i[t]$ is the transmit power of device $i$ during time interval $t$. Therefore, the rate of IoT device $i$ with respect to UAV $u$ during time slot $t$ can be written as ${R}_{iu}[t]=B_{iu}[t]\log _{2}\left(1+\Gamma_{iu}[t]\right),$ where $B_{iu}[t]$ is the allocated bandwidth between device $i$ and UAV $u$ during time slot $t$.




 During its flight, a UAV makes stops to collect data from subsets of IoT devices. To ensure efficient and rapid data transmission between device $i$ and UAV $u$, the data rate denoted as $R_{iu}$ for the device-UAV pair should exceed a predefined threshold value $R^{\text{min}}$. This threshold, $R^{\text{min}}$, is set high enough to enable nearly instantaneous transmission of all accumulated data updates. Consequently, the data collection time is considered negligible compared to the flight time. For the purpose of our analysis, we assume that UAVs maintain a constant flight speed denoted as $V$. Hence, the total flight time  for UAV $u$, denoted by $\zeta_u$, can be expressed as $ \zeta_{u}(\!\boldsymbol{x}_u,\boldsymbol{y}_u\!)\!=\!\!\!
   \sum_{t=1}^{K}\!\!\frac{\sqrt{(x_{u}[t\!+\!1]\!-\!x_{u}[t])^2\!+\!(y_{u}[t\!+\!1]\!-\!y_{u}[t])^2}}{V}.$   

Our objective is to optimize the trajectories of UAVs and the subsets of visited IoT devices in order to minimize the AoI across all devices subject to some quality of service constraints. To achieve this goal, we will first present the mathematical formulation of the problem. Subsequently, we will propose an efficient distributed approach to solve this optimization problem.

\section{Problem Formulation }\label{Prob}
The main aim of our paper is to minimize the AoI across all devices throughout the time span $T$. To achieve this objective, we focus on optimizing the stopping locations of UAVs over time and the selection of devices for data collection. The optimization problem is formulated as follows.
 
\begin{mini!}
{ \boldsymbol{\alpha},\boldsymbol{x},\boldsymbol{y}}{\sum\limits_{i \in \mathcal{I}} f_i(\boldsymbol{\alpha}_i)
,\label{objective}}
{\label{GeneralOptimizati}}{}
\addConstraint{R_{iu}[t]\geq \alpha_{iu}[t]R^{\rm min},\;  \forall u \in \mathcal{U}, \forall i \in \mathcal{I} \label{Rate}}{}
\addConstraint{}{  \forall t \in \{0,\dots, K-1\} }\nonumber
{}{}
\addConstraint{ 
\sum \limits_{u\in \mathcal{U}}\alpha_{iu}[t]\leq 1, \; \forall i \in \mathcal{I},  \forall t \in \{0,\dots, K-1\}
\label{Association2}}
{}{}
\addConstraint{ 
\zeta_u(\boldsymbol{x}_u,\boldsymbol{y}_u)\leq \zeta_u^{\rm max}, \;  \forall  \in \mathcal{U}
\label{Time}}
{}{}
\addConstraint{ 
0\leq x_u[t]\leq x^{\rm max}, \;  \forall u \in \mathcal{U}, t \in \{0,K-1\}
\label{xposition}}
{}{}
\addConstraint{ 
0\leq y_u[t]\leq y^{\rm max}, \;  \forall u \in \mathcal{U}, t \in \{0,K-1\}
\label{yposition}}
{}{}
\addConstraint{ 
\alpha_{iu}[t]\in \{0,1\}, \;  \forall i \in \mathcal{I}, u \in \mathcal{U}, t \in \{0,K-1\}
\label{alpha}}
{}{}
\end{mini!}

Constraint~(\ref{Rate}) ensures that the data rate between each UAV and its associated IoT device remains above the predefined threshold $R^{\text{min}}$. Constraint~(\ref{Association2}) ensures that each IoT device transmits data to only one UAV at a time. Constraint~(\ref{Time}) guarantees that each UAV's flight time, $\zeta_u$, does not exceed its maximum allowed flight time, $\zeta_u^{\text{max}}$, which is determined by its energy budget. Constraints~(\ref{xposition}) and (\ref{yposition}) restrict the movement of UAVs within the studied area. Constraint~(\ref{alpha}) reflects the binary nature of the association variables.

The studied optimization problem involves a combination of binary and continuous variables, as well as nonlinear objective function and constraints. As a result, it falls under the category of Mixed-Integer-Nonlinear-Programming (MINLP) problems, which are known to be challenging to solve. Traditional MINLP optimization algorithms such as Branch-and-Bound are not appropriate to solve such problems due to their exponential convergence time and the need for global knowledge of the environment and its dynamics. Additionally, classical ML algorithms, including standard RL, struggle to scale with the large dimension of the variable set and are not well-suited for multi-agent environments. To address these challenges, we leverage the MAPPO  framework to efficiently solve the optimization problem and find optimal trajectories and device subsets for data collection by the UAVs.

\section{A Decentralized MAPPO Implementation for AoI Minimization}\label{algo}
The optimization of trajectories and subsets of served IoT devices in a multi-UAV system is a complex task due to the interdependence of UAV policies and their shared objective of minimizing the global AoI. Each UAV's decision to visit and collect data from a device depends on whether other UAVs have previously collected data from that device or not. However, achieving a fully synchronized system where UAVs possess complete global knowledge of the environment, including the positions and historical actions of other UAVs, is challenging in practical scenarios. 


To cope with the challenges of optimizing a common objective in a multi-UAV system, we employ the MAPPO approach, which extends the PPO algorithm. The MAPPO framework adopts a multi-agent learning perspective, allowing individual UAVs to make decisions independently.  The MAPPO approach utilizes an actor-critic architecture. The actor is a neural network that approximates the agent's policy and is responsible for decision-making, while the critic is another neural network that evaluates these decisions by estimating a value function.

To efficiently train the MAPPO agents, we utilize the DTDE scheme. Unlike CTDE where all the agents policies are fed into one centralized critic network, the DTDE setup assumes a critic architecture for each agent. Specifically, each critic learns an individual value function that considers the local actions of its respective agent. On the other hand, the actor network is employed by each UAV to calculate its policy, relying solely on its own local observations. 

 To achieve MARL based on the DTDE scheme, the objective function needs to be decomposed so that each agent can compute its own reward. However, directly decomposing the original objective in (\ref{objective}) is challenging due to the need for information sharing between agents about which IoT devices have been visited previously and their data have been collected. 
 
 To overcome this limitation, we design a new objective function in (\ref{objective2}) that can easily be decomposed into individual rewards. The proposed function involves the number of served devices at time step $t$ which is weighted by their respective AoI at the previous time step, $t-1$. This modification ensures that agents take into account the cumulated age of data when making decisions, leading to a more effective reduction in the overall AoI.


Consequently, an alternative optimization problem, which can be seen as an approximation of the original problem (\ref{GeneralOptimizati}), is formulated as follows

\begin{maxi!}
{ \boldsymbol{\alpha},\boldsymbol{x},\boldsymbol{y}}{\sum \limits_{u \in \mathcal{U}}
   \sum \limits_{t=1}^K\sum\limits_{i \in \mathcal{I}}\alpha_{iu}[t]\sum\limits_{n=0}^{\lfloor \frac{K}{k_i} \rfloor}w^n[t]A_{i}^n[t-1]
,\label{objective2}}
{\label{ApproxOptimizati}}{}
\addConstraint{}{ (\ref{Rate}), (\ref{Association2}), (\ref{Time}), (\ref{xposition}), (\ref{yposition}) \text{ and } (\ref{alpha}). }
{}{}
\end{maxi!}

It can be seen from this formulation that the new objective function is composed of the sum of individual rewards of each UAV, where the target of each UAV is to maximize the number of visited devices weighted by their AoI in the previous time step. Specifically, each UAV at time $t$, will seek to visit devices with the highest AoI at time $t-1$, in order to reset their AoI to $0$ at time $t$.
This modified problem serves as a basis for developing an efficient and distributed MARL approach using the MAPPO framework.

The proposed decentralized MAPPO algorithm, referred to as MAPPO-AoI-Dec, involves the following elements.
\begin{itemize}
    \item \textbf{Agents}: the set of UAVs.

    \item \textbf{Actions}: 
    At each time step, two decisions need to be made by each UAV. The first decision involves moving along either the $x$ or $y$ axis, with actions  \textit{up, down, right, left}, or \textit{staying in the same place}. The second decision involves selecting a subset of devices from which data is collected. Hence, an action of a UAV $u$ at time $t$ can be represented as a tuple $\boldsymbol{a}_u[t] \in \{0,1,2,3,4\}\times \{0,1\}^I$. In this representation, the first component, which can take values from ${0,1,2,3,4}$, indicates whether the UAV stays at the same position ($0$) or moves in a specific direction: up ($1$), down ($2$), right ($3$), or left ($4$). The second component is a binary vector of length $I$ that specifies the devices selected for data collection.

    
    \item \textbf{States}: The state of UAV $u$ at time step $t$ is represented by its location, denoted as $s_u[t]$. To encompass all possible states, the studied area is divided into a 2D grid with a granularity of $\delta$, denoted as $\mathcal{G}_\delta$. Therefore, the set of states for a UAV is given by $\mathcal{S}=\mathcal{G}_\delta$. The global state vector of all UAVs at time step $t$ is denoted as $\boldsymbol{s}[t]=(s_1[t],\dots,s_U[t])$. 
    
    \item  \textbf{Rewards}: The reward of a UAV $u$ at time $t$, $r_{u}[t]$, is given by the sum of its visited devices weighted by their AoI at the previous time interval, i.e.,
    $r_{u}[t]=\sum\limits_{i \in \mathcal{I}}\alpha_{iu}[t]\sum\limits_{n=0}^{\lfloor \frac{K}{k_i} \rfloor}w^n[t]A_{i}^n[t-1].$
    
    \item \textbf{Policy function}: We define the policy function $\boldsymbol{\pi}_{\boldsymbol{\theta}_{u}}(\boldsymbol{a}_{u},s_{u})$ as the actor network that takes as input the local state $s_u$ of UAV $u$ and outputs a probability distribution over its possible actions $\boldsymbol{a}_{u}$. The policy function is parameterized by the vector $\boldsymbol{\theta}_{u}$ and is used to generate the UAV's strategy based on its own local state. The UAV then selects an action according to this probability distribution.

    \item \textbf{Value function}: The value function $V_{\boldsymbol{\phi}}(\boldsymbol{s}_u)$ is represented by the critic network and is parameterized by the vector $\boldsymbol{\phi}_u$. The critic network is trained to estimate the expected future reward that can be achieved by a UAV $u$ starting from a given state $\boldsymbol{s}_u$. The critic network takes as input the state $\boldsymbol{s}_u$ and provides a scalar output representing the estimated reward that the UAV can expect to receive in the future if it follows the policy $\boldsymbol{\pi}_{\boldsymbol{\theta}_{u}}$. The goal of the critic network is to find the optimal strategy that maximizes the expected reward of UAV $u$.
\end{itemize}
The MAPPO-AoI-Dec algorithm is structured similarly to the PPO algorithm used in single-agent settings~\cite{Schulman2017ProximalPO}. In this approach, each UAV's actor network updates its policy parameters using the PPO algorithm based on its own local state and reward information. At each training step $t$, when presented with its current state $s_u[t]$, the actor network computes an action $\boldsymbol{a}_u[t]$ according to its policy $\boldsymbol{\pi}_{\boldsymbol{\theta_u}}(\boldsymbol{a}_u[t], s_u[t])$. This selected action is then executed, leading to a new state $s_u[t+1]$ and a corresponding local reward $r_u[t]$.

During the training process, the MAPPO-AoI-Dec algorithm individually trains critic networks to estimate the value function for each UAV agent $u$. Each critic network takes only the local state $s_u$ of its corresponding agent as input and outputs the estimated value $V_{\boldsymbol{\phi_u}}(s_u)$. The training procedure of MAPPO-AoI-Dec involves alternating between updating the actor and critic networks until convergence, as illustrated in Fig.~\ref{fig:dtde}. Each UAV trains its actor network using the PPO approach, minimizing a loss function based on its local state $s_u[t]$ at each training step $t$. The actor network samples an action $\boldsymbol{a}_u[t]$ from its policy $\boldsymbol{\pi}_\theta(\boldsymbol{a}_u[t],s_u[t])$, leading to a new state $s_u[t+1]$. Subsequently, the local state vector of each UAV is fed into each of their respective critic networks, which are also trained using the PPO approach. Moreover, before moving to the next state and observing the reward, the MAPPO-AoI-Dec algorithm ensures that the constraints of the optimization problem are satisfied.




\section{Simulation Results}\label{simu}
\begin{table}[ht]
\centering
\begin{tabular}{|l|l||l|l|}
\hline
\textbf{Parameter} & \textbf{Value} & \textbf{Parameter} & \textbf{Value}   \\ \hline
$I$       & $25$      &$U$       & $3$  \\ \hline
$x_{max}$       & $1000m$ & $y_{max}$       & $1000m$  \\ \hline
$H_u$       & $[80,100]m$   & $R_{min}$       & $150$ Kbit/s  \\ \hline
$\tau$       & $3ms$   &$k_i$       & $[1,5]$  \\ \hline
$B_{i,u}$       & $[1.5,2]$ GHz   &$\sigma^{2}$       & $-120dBm$  \\ \hline
$P_{i}$       & $[0,1]$ mW   &$T$       & $500\tau$  \\ \hline
\end{tabular}
\caption{Simulation setup}
\label{tab:my-table-parameter}
\end{table}

\begin{figure*}[ht]
\centering
\minipage{0.31\textwidth}
  \includegraphics[width=\textwidth]{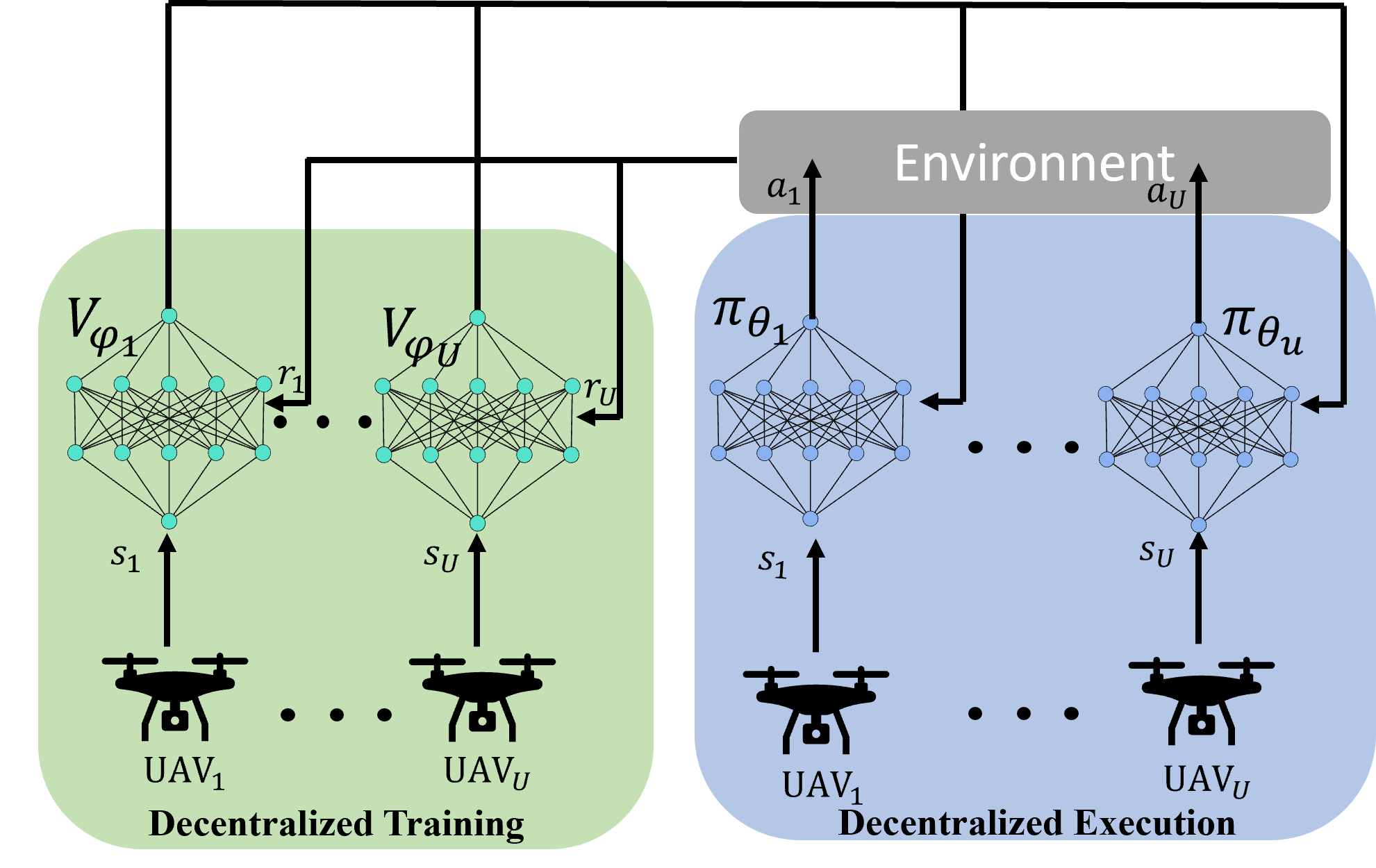}
   \caption{An illustration of a DTDE approach.}
    \label{fig:dtde}
\endminipage \hfill
\minipage{0.31\textwidth}
  \includegraphics[width=\textwidth]{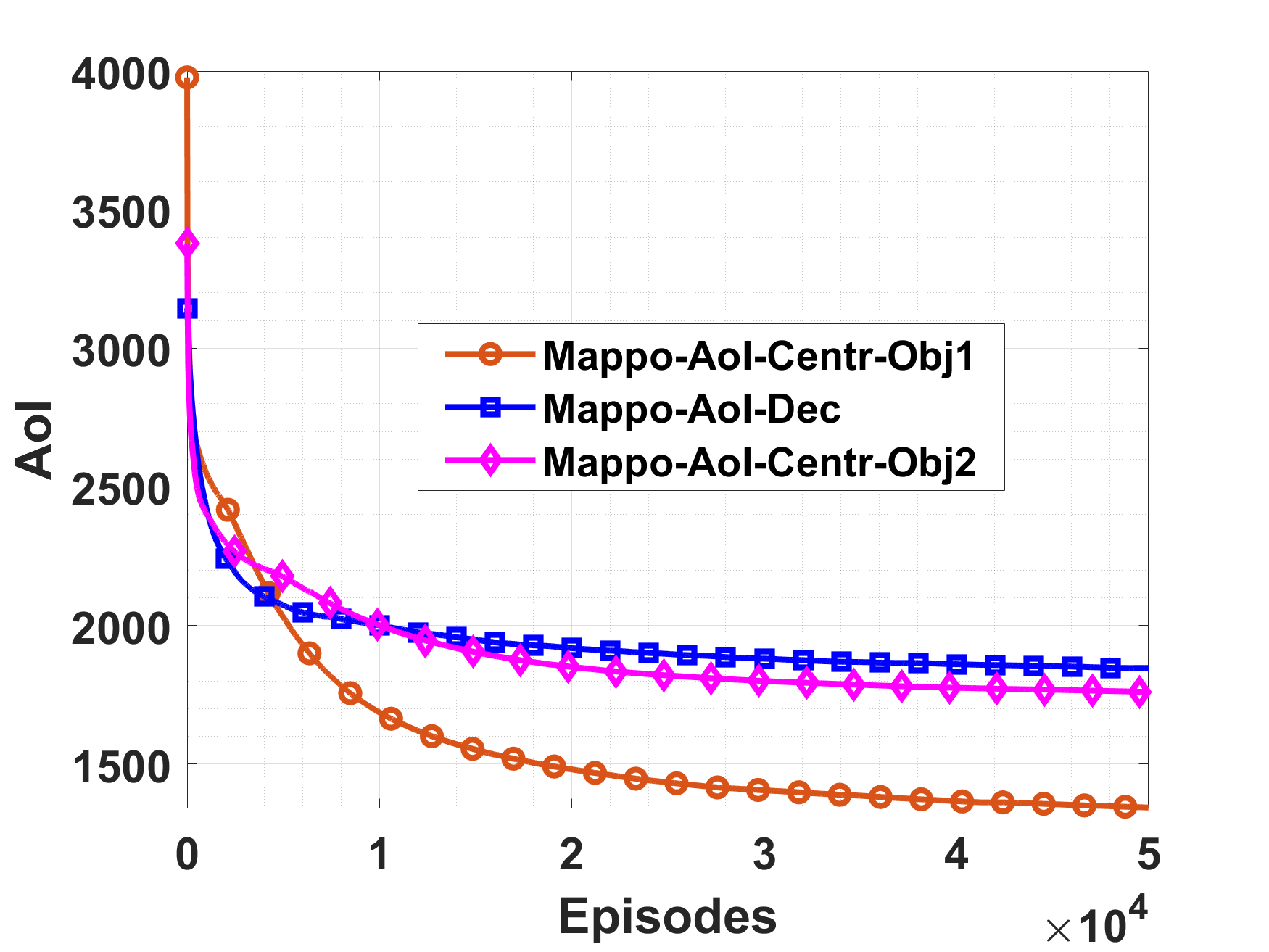}
\caption{AoI over \\episodes.}
 \label{fig:AoI}
\endminipage\hfill
\minipage{0.31\textwidth}
  \includegraphics[width=\textwidth]{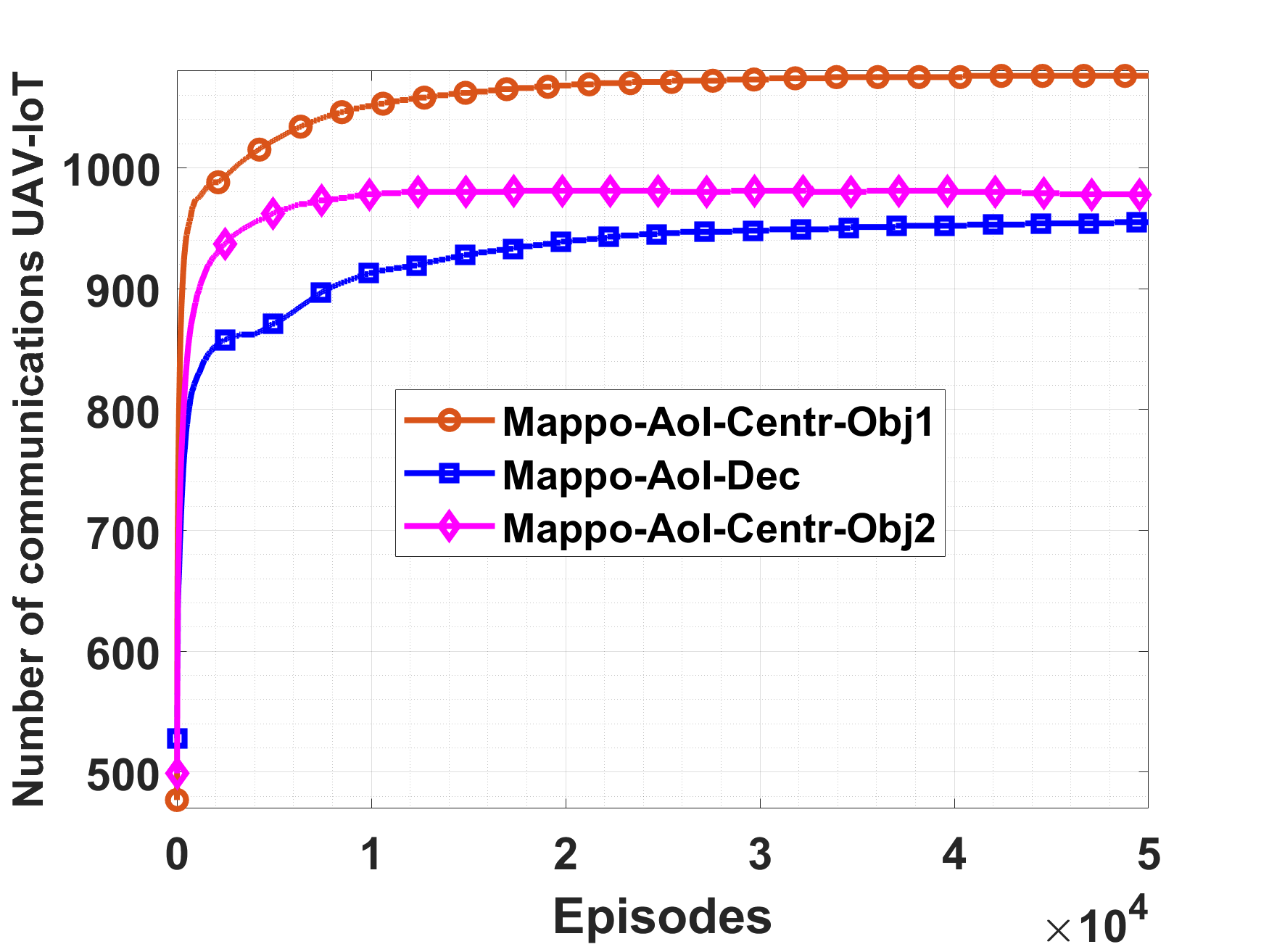}
   \caption{ Total number of UAV-IoT communications over time.}
    \label{fig:iotnondistinct}
\endminipage 
\end{figure*}

We evaluate the performance of our proposed approach through simulations with a $1000m\times1000m$ area, $25$ IoT devices randomly deployed and $3$ UAVs for data collection. We assume that the UAVs fly at altitudes $80m$, $90m$ and $100m$, and maintain a constant speed of $15m/s$. The mission time is set to $\zeta_u^{\rm max}=500\tau$, where $\tau=3ms$. The data generation period of IoT devices varies randomly between $1\tau$ and $5\tau$. Each device is assigned a fixed bandwidth randomly selected from the range $[1.5,2]$ GHz and a constant power between $[0,1]$ mWatt. To meet the quality of service constraint, a minimum rate of $150$ Kbit/s is ensured. The simulation parameters are summarized in Table \ref{tab:my-table-parameter}. Furthermore, we set $w^n[t]=\gamma^{t-n}$ where $\gamma=0.8$.

In our simulations, we compare two approaches for training MARL agents in the context of UAV-assisted data collection.

\begin{itemize}
    \item {\textbf{MAPPO-AoI-Dec}} is the proposed decentralized approach which solves problem~(\ref{objective2}) with a decentralized implementation of MAPPO.
    \item \textbf{MAPPO-AoI-Centr-Obj1} is a centralized MAPPO similar to the one in~\cite{naby2023muti} that solves the original optimization problem.
    \item \textbf{MAPPO-AoI-Centr-Obj2} is a centralized MAPPO as in~\cite{naby2023muti} that solves the approximate problem~(\ref{objective2}) using a centralized training.
\end{itemize}

\begin{figure*}[!h]
\minipage{0.31\textwidth}
  \includegraphics[width=\textwidth]{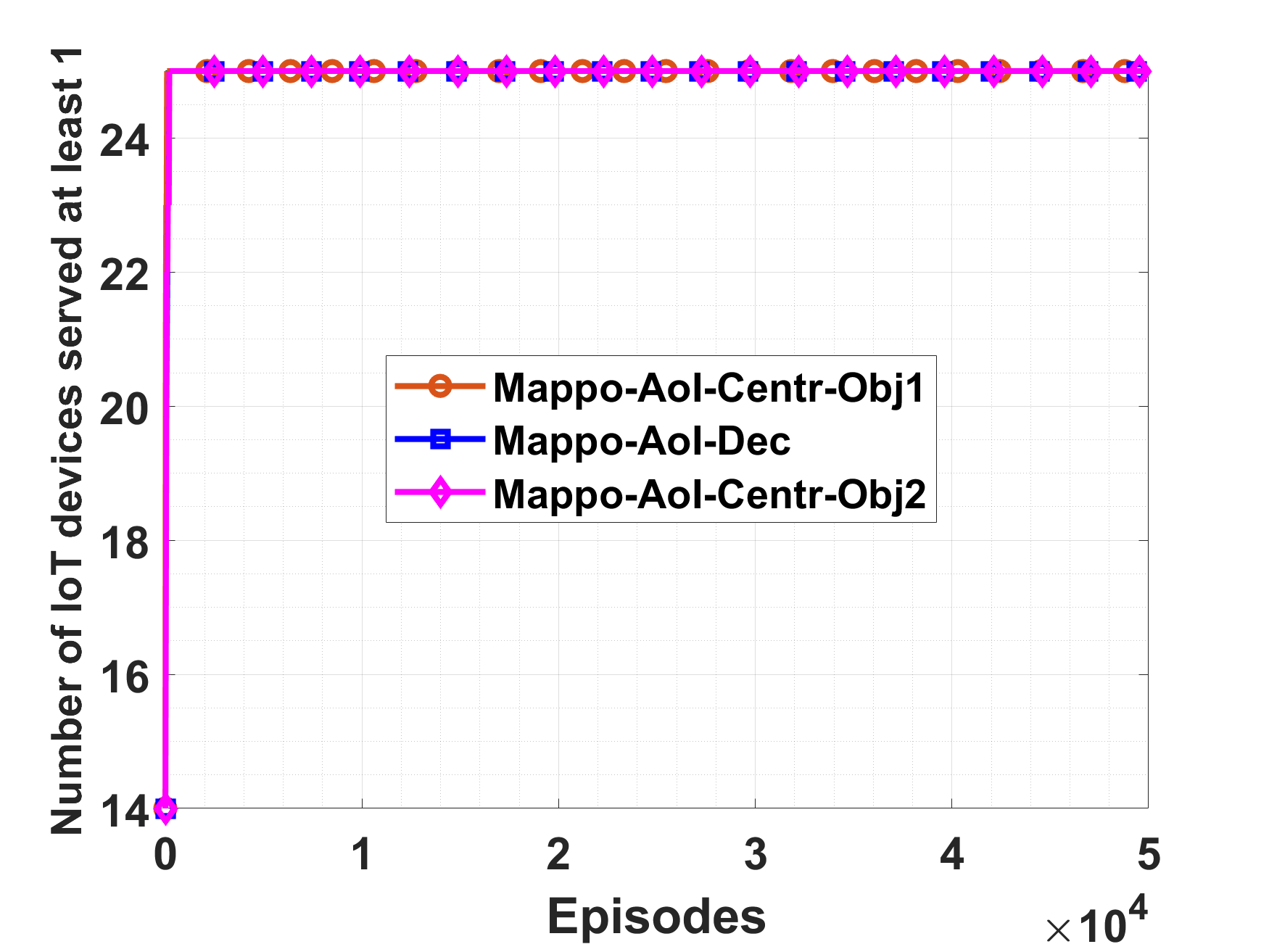}
   \caption{Number of served IoT devices at least once over episodes.}
    \label{fig:iotdistinct}
\endminipage\hfill
\minipage{0.31\textwidth}
  \includegraphics[width=\textwidth]{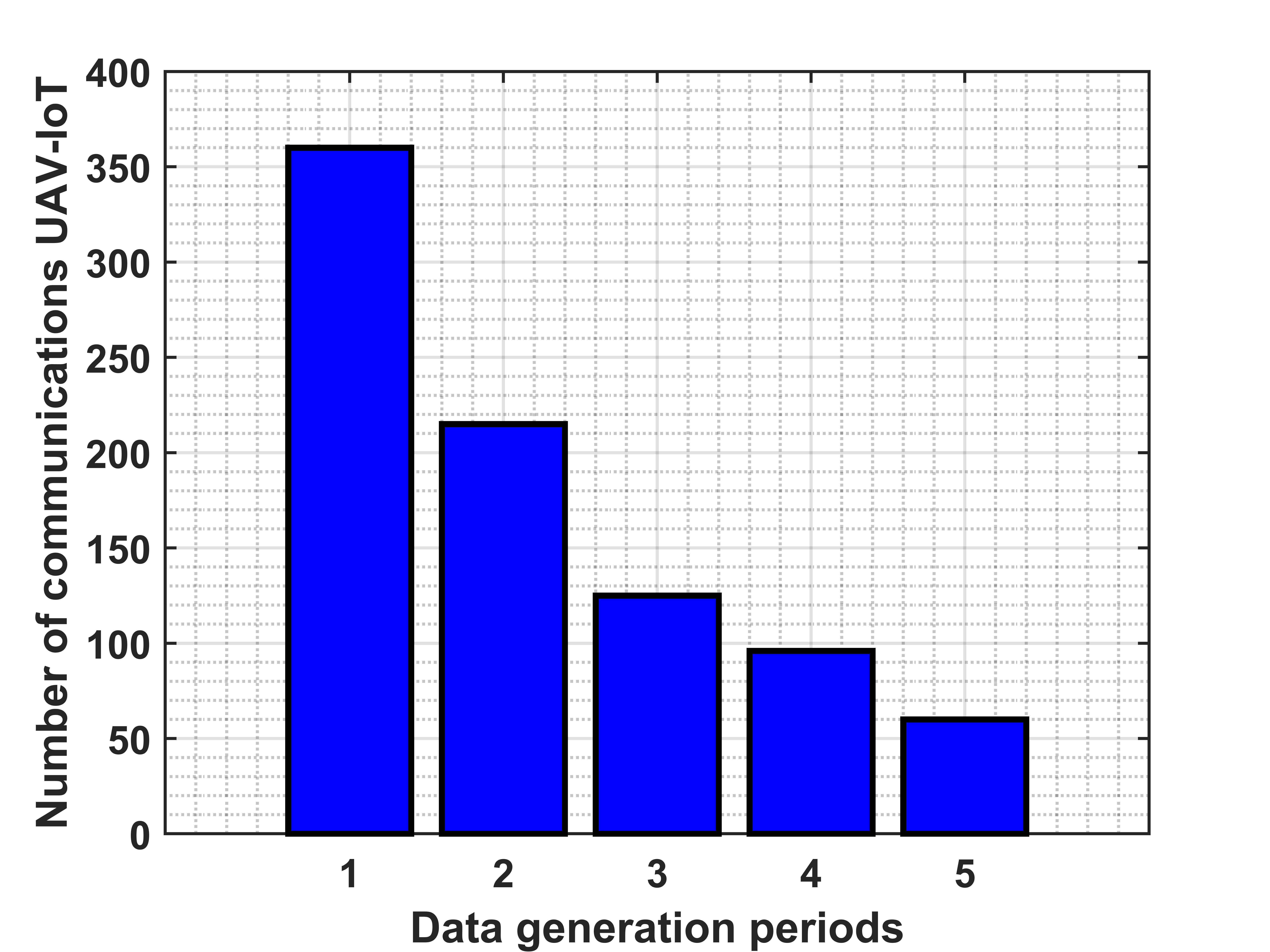}
   \caption{MAPPO-AoI-Dec with the approximate problem.}
   \label{fig:gendata}
\endminipage\hfill
\minipage{0.31\textwidth}
  \includegraphics[width=\textwidth]{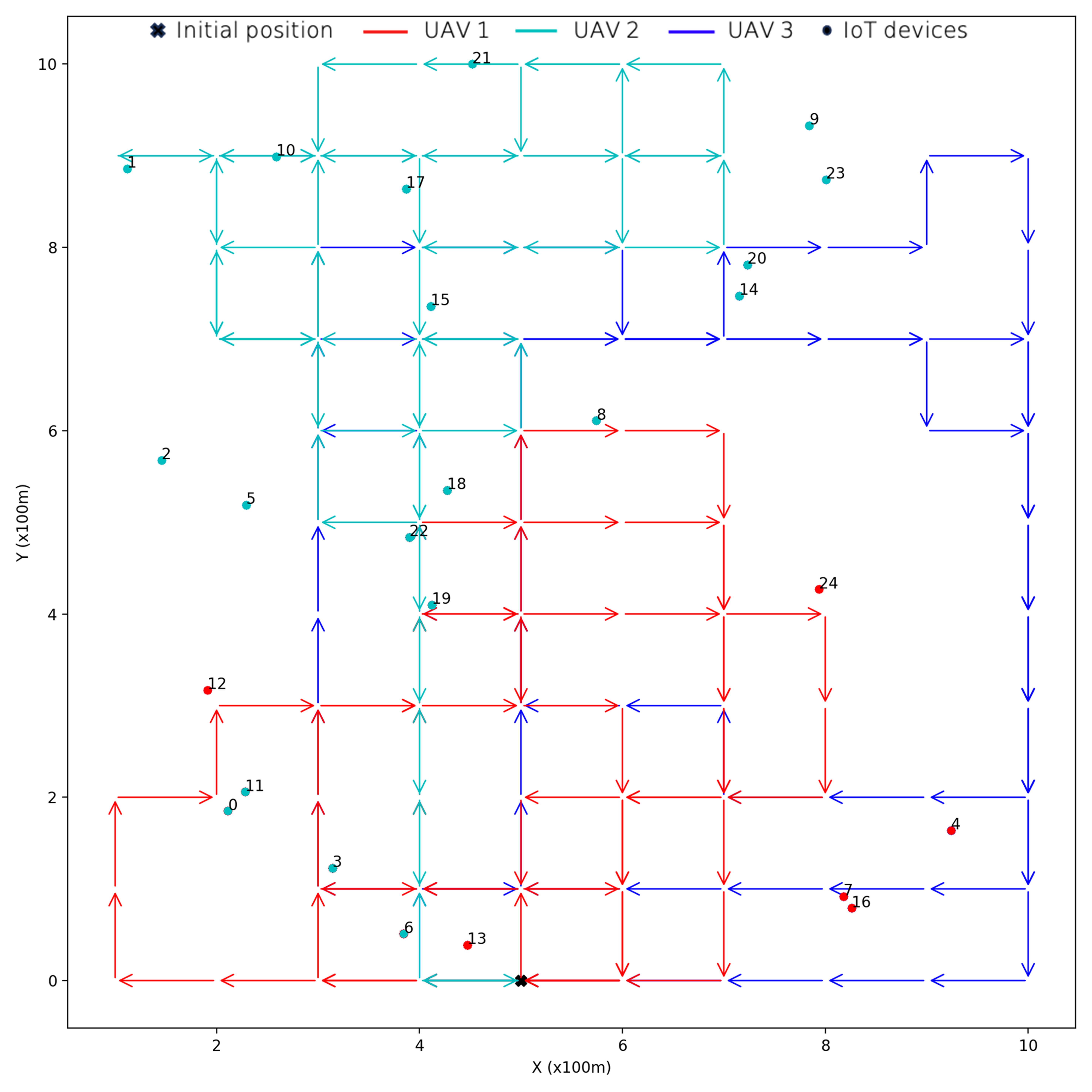}
   \caption{UAVs trajectories.}
   \label{fig:traj}
\endminipage
\end{figure*}


In Fig. \ref{fig:AoI}, we can observe the total AoI of the system throughout the training episodes. It can be seen from the plot that all three studied schemes achieve nearly the same AoI at the end of the training. However, as expected, the MAPPO-AoI-Centr-Obj1 scheme demonstrates slightly better performance due to its optimization of the exact target function depicted in the figure. It is followed by MAPPO-AoI-Centr-Obj2 which also benefits from the global knowledge of the network, contributing to its improved results. Finally, MAPPO-AoI-Dec achieves promising results although the agents, for this scheme, do not share any information with each other.



In Fig. \ref{fig:iotnondistinct}, we plot the total number of communications between IoT devices and UAVs during their mission. This number represents the total number of times an IoT device sends its collected data to a UAV. It can be noticed that the number of data transmissions reaches its highest values for the centralized schemes, aligning with the results shown in Fig. \ref{fig:AoI}. 

In Fig.~\ref{fig:iotdistinct}, we plot the number of IoT devices visited at least once. As depicted, all devices have had their data collected at least once for the three schemes. This is expected since leaving a device unvisited would significantly increase the total AoI, making it essential to collect data from all devices at least once.

Fig. \ref{fig:gendata} illustrates the relationship between the number of UAV-IoT communications and the data generation period in the proposed scheme. The figure highlights that IoT devices with shorter data generation periods experience a greater frequency of interactions with UAVs, making them more frequently visited.
Finally, the trajectories of all UAVs are depicted in Fig. \ref{fig:traj}. Commencing their flight from a designated point such as CDS, these UAVs navigate to visit all IoT devices to maximize the freshness of the collected data.

\section{Conclusion}\label{Conc}
Our paper focused on the deployment of UAVs for time-sensitive data collection and quantifying data freshness using the AoI metric. We formulated the problem as a MINLP and proposed a MAPPO approach to efficiently solve it. To ensure scalability and effectiveness, we developed a decentralized method for training MARL agents and compared it to a centralized approach. Simulation results showed that the decentralized approach achieved promising results. 

\section*{Acknowledgment}
This document has been produced with the financial assistance of the European Union (Grant no. DCI-PANAF/2020/420-028), through the African Research Initiative for Scientific Excellence (ARISE), pilot programme. ARISE is implemented by the African Academy of Sciences with support from the European Commission and the African Union Commission.

\balance

\bibliographystyle{IEEEtran}
\bibliography{WCNC}

\begin{thebibliography}{10}
\providecommand{\url}[1]{#1}
\csname url@samestyle\endcsname
\providecommand{\newblock}{\relax}
\providecommand{\bibinfo}[2]{#2}
\providecommand{\BIBentrySTDinterwordspacing}{\spaceskip=0pt\relax}
\providecommand{\BIBentryALTinterwordstretchfactor}{4}
\providecommand{\BIBentryALTinterwordspacing}{\spaceskip=\fontdimen2\font plus
\BIBentryALTinterwordstretchfactor\fontdimen3\font minus \fontdimen4\font\relax}
\providecommand{\BIBforeignlanguage}[2]{{%
\expandafter\ifx\csname l@#1\endcsname\relax
\typeout{** WARNING: IEEEtran.bst: No hyphenation pattern has been}%
\typeout{** loaded for the language `#1'. Using the pattern for}%
\typeout{** the default language instead.}%
\else
\language=\csname l@#1\endcsname
\fi
#2}}
\providecommand{\BIBdecl}{\relax}
\BIBdecl

\bibitem{bajracharya20226g}
R.~Bajracharya, R.~Shrestha, S.~Kim, and H.~Jung, ``{6G} {NR-U} based wireless infrastructure {UAV}: Standardization, opportunities, challenges and future scopes,'' \emph{IEEE Access}, vol.~10, pp. 30\,536--30\,555, 2022.

\bibitem{hammouti2018air}
H.~E. Hammouti and M.~Ghogho, ``Air-to-ground channel modeling for uav communications using 3d building footprints,'' in \emph{Ubiquitous Networking: 4th International Symposium, UNet 2018, Hammamet, Tunisia, May 2--5, 2018, Revised Selected Papers 4}.\hskip 1em plus 0.5em minus 0.4em\relax Springer, 2018, pp. 372--383.

\bibitem{kosta2017age}
A.~Kosta, N.~Pappas, V.~Angelakis \emph{et~al.}, ``Age of information: A new concept, metric, and tool,'' \emph{Foundations and Trends{\textregistered} in Networking}, vol.~12, no.~3, pp. 162--259, 2017.

\bibitem{10333749}
M.~N. Ndiaye, E.~H. Bergou, and H.~E. Hammouti, ``Ensemble dnn for age-of-information minimization in uav-assisted networks,'' in \emph{2023 IEEE 98th Vehicular Technology Conference (VTC2023-Fall)}, 2023, pp. 1--6.

\bibitem{arulkumaran2017deep}
K.~Arulkumaran, M.~P. Deisenroth, M.~Brundage, and A.~A. Bharath, ``Deep reinforcement learning: A brief survey,'' \emph{IEEE Signal Processing Magazine}, vol.~34, no.~6, pp. 26--38, 2017.

\bibitem{Cui2020MultiAgentRL}
J.~Cui, Y.~Liu, and A.~Nallanathan, ``Multi-agent reinforcement learning-based resource allocation for {UAV} networks,'' \emph{IEEE Transactions on Wireless Communications}, vol.~19, pp. 729--743, 2020.

\bibitem{Seid2021MultiAgentDF}
A.~M. Seid, G.~O. Boateng, B.~Mareri, G.~Sun, and W.~Jiang, ``Multi-agent drl for task offloading and resource allocation in multi-{UAV} enabled {IoT} edge network,'' \emph{IEEE Transactions on Network and Service Management}, vol.~18, pp. 4531--4547, 2021.

\bibitem{9900429}
S.~Araf, A.~S. Saha, S.~H. Kazi, N.~H. Tran, and M.~G.~R. Alam, ``{UAV} assisted cooperative caching on network edge using multi-agent actor-critic reinforcement learning,'' \emph{IEEE Transactions on Vehicular Technology}, pp. 1--16, 2022.

\bibitem{Lin2021DecentralizedPD}
J.~Lin, H.-T. Chiu, and R.-H. Gau, ``Decentralized planning-assisted deep reinforcement learning for collision and obstacle avoidance in {UAV} networks,'' \emph{2021 IEEE 93rd Vehicular Technology Conference (VTC2021-Spring)}, pp. 1--7, 2021.

\bibitem{Hu2021DistributedMM}
Y.~Hu, M.~Chen, W.~Saad, H.~V. Poor, and S.~Cui, ``Distributed multi-agent meta learning for trajectory design in wireless drone networks,'' \emph{IEEE Journal on Selected Areas in Communications}, vol.~39, pp. 3177--3192, 2021.

\bibitem{Hu2022DistributedAD}
Y.~Hu, X.~Wang, and W.~Saad, ``Distributed and distribution-robust meta reinforcement learning (d2-rmrl) for data pre-storing and routing in cube satellite networks,'' \emph{ArXiv}, vol. abs/2206.06568, 2022.

\bibitem{Schulman2017ProximalPO}
J.~Schulman, F.~Wolski, P.~Dhariwal, A.~Radford, and O.~Klimov, ``Proximal policy optimization algorithms,'' \emph{ArXiv}, vol. abs/1707.06347, 2017.

\bibitem{Yu2021TheSE}
C.~Yu, A.~Velu, E.~Vinitsky, Y.~Wang, A.~M. Bayen, and Y.~Wu, ``The surprising effectiveness of mappo in cooperative, multi-agent games,'' \emph{ArXiv}, vol. abs/2103.01955, 2021.

\bibitem{naby2023muti}
M.~N. Ndiaye, E.~Bergou, and H.~E. Hammouti, ``Muti-agent proximal policy optimization for data freshness in {UAV}-assisted networks,'' in \emph{IEEE Communication Conference workshops}, 2023.

\bibitem{ndiaye2022age}
M.~N. Ndiaye, E.~Bergou, M.~Ghogho, and H.~E. Hammouti, ``Age-of-updates optimization for {UAV}-assisted networks,'' in \emph{IEEE Global Communications Conference}, 2022, pp. 450--455.

\bibitem{AbdElmagid2019DeepRL}
M.~A. Abd-Elmagid, A.~Ferdowsi, H.~S. Dhillon, and W.~Saad, ``Deep reinforcement learning for minimizing age-of-information in {UAV}-assisted networks,'' \emph{2019 IEEE Global Communications Conference (GLOBECOM)}, pp. 1--6, 2019.

\bibitem{kadota2021age}
I.~Kadota and E.~Modiano, ``Age of information in random access networks with stochastic arrivals,'' in \emph{IEEE INFOCOM 2021-IEEE Conference on Computer Communications}.\hskip 1em plus 0.5em minus 0.4em\relax IEEE, 2021, pp. 1--10.

\bibitem{bedewy2019age}
A.~M. Bedewy, Y.~Sun, and N.~B. Shroff, ``The age of information in multihop networks,'' \emph{IEEE/ACM Transactions on Networking}, vol.~27, no.~3, pp. 1248--1257, 2019.

\end{thebibliography}

\end{document}